\documentclass[12pt]{amsart}
\usepackage[all]{xy}

\newtheorem{thm}{Theorem}[section]
\newtheorem{prob}[thm]{Problem}

\theoremstyle{definition}
\newtheorem{defn}[thm]{Definition}

\theoremstyle{remark}
\newtheorem{rem}[thm]{Remark}

\renewcommand{\inf}{P_\infty(\N)}

\renewcommand{\d}{{\mathfrak d}}
\newcommand{\p}{{\mathfrak p}}

\newcommand{\CH}{the Continuum Hypothesis}

\newcommand{\w}{\omega}

\renewcommand{\t}{\mathfrak{t}}

\newcommand{\bq}{\begin{quote}}
\newcommand{\eq}{\end{quote}}

\newcommand{\sone}{\mathsf{S}_1}

\newcommand{\cF}{\mathcal{F}}
\newcommand{\cG}{\mathcal{G}}

\newcommand{\psin}{pseudo-intersection}

\newcommand{\naturals}{{\mathbb N}}
\newcommand{\N}{\naturals}
\newcommand{\as}{\subseteq^*}
\newcommand{\sm}{\setminus}

\newcommand{\by}[2]{\par\hfill\emph{#1}, #2}
\newcommand{\Tau}{\mathrm{T}}
\newcommand{\CE}{\textsc{CE}}

\newcommand{\be}{\begin{enumerate}}
\newcommand{\ee}{\end{enumerate}}
\newcommand{\bi}{\begin{itemize}}
\newcommand{\ei}{\end{itemize}}
\renewcommand{\i}{\item}

\newcommand{\arx}[1]{\texttt{http://arxiv.org/abs/#1}}

\title[$\mathcal{SPM}$ Bulletin \textbf{5} (May 2003)]{%
$\mathcal{SPM}$ Bulletin\\[0.5cm]
Issue number 5: May 2003 \CE{}}

\begin{document}
\maketitle

\tableofcontents

\begin{rem}
(The current remark is added several months after this issue was released.)
We were informed by Shelah that he found a gap in his solution to the problem
discussed below.
\emph{The minimal tower problem is still open.}
The positive part of this fact is that  
this problem will continue to inspire mathematicians 
and yield related results. We leave the rest of this issue unchanged 
for documentational reasons.
\end{rem}

\section{Editor's note}

This is a special issue dedicated to the announcement
of Shelah's recent solution of the Minimal Tower problem,
one of the oldest and most important problems in infinite
combinatorics which also motivated some new studies in
SPM (see fourth issue of this bulletin).
We give some background and personal perspectives on the 
problem and its solution.

We decided not to include additional research 
announcements in this issue, so to let them draw
the attention they deserve in the coming issue.

The first issues of this bulletin,
which contain general information (first issue),
basic definitions, research announcements, and open problems (all issues) are available online:
\be
\i First issue: \arx{math.GN/0301011}
\i Second issue: \arx{math.GN/0302062}
\i Third issue: \arx{math.GN/0303057}
\i Fourth issue: \arx{math.GN/0304087}
\ee

{\small
\subsection*{Contributions}
Please submit your contributions (announcements, discussions, and open problems)
by e-mailing us. It is preferred to write them
in \LaTeX{} (or otherwise in some other \TeX{} format
or plain text).
The authors are urged to use as standard notation as possible, or otherwise give
a reference to where the notation is explained.
Contributions to this bulletin would not require any transfer of copyright,
and material presented here can be published elsewhere.

\subsection*{Subscription}
To receive this bulletin (free) to your
e-mailbox, e-mail us.
}

\medskip

\by{Boaz Tsaban}{tsaban@math.huji.ac.il}

\hfill \texttt{http://www.cs.biu.ac.il/\~{}tsaban}

\section{The Minimal Tower problem solved: A personal perspective}
Yesterday (May 22) I met Saharon Shelah in the Hebrew
University, and he told me the astonishing news
that he has just solved the Minimal Tower problem!
The details of the proof are currently being checked, 
but he believes that the proof is correct and the method
of proof is closely related to a method he used to settle
another major open problem two years ago, so we assume in
the sequel that the problem is indeed settled.
Below are some details on the problem and \emph{personal}
(mathematical and nonmathematical) impressions on the 
history and motivation of Shelah's working on the problem.

This problem asks, for two cardinal characteristics of the
continuum called $\p$ and $\t$ (to be defined below),
whether it is provable that $\p=\t$.
Shelah proved that the answer is ``No'': 
It is consistent with the usual axioms of mathematics
(ZFC) that $\p<\t$.

A simple formulation of the \emph{Minimal Tower problem (MTP)}
is as follows:
Write $A\as B$ for ``$A\sm B$ is finite''.
A set $A$ is a \emph{pseudo-intersection} of
a family $\cF$ of sets if for each $B\in\cF$,
$A\as B$.
\begin{prob}[Solved by Shelah, May 2003]\label{MTP}
Does $(*)_\kappa$ hold for each cardinal $\kappa$?
\bi
\i[$(*)_\kappa$] Assume that each $\as$-chain $\cF$ of infinite sets of natural numbers
with $|\cF|\le\kappa$ has a \psin{}. 
Then for each family $\cG$ of infinite sets of natural numbers
such that $|\cG|\le\kappa$ and $\cG$ 
is closed under taking finite intersections, $\cG$ has a \psin{}.
\ei
\end{prob}
It is not difficult to see that
for $\kappa\le\aleph_0$ the conclusion in the assertion $(*)_\kappa$ 
is true and for $\kappa=2^{\aleph_0}$ the assumption in $(*)_\kappa$ is false.
Thus, if \CH{} holds, then the answer is positive.
But it is not clear what happens if we do not assume \CH{} 
(or any other similar assumption).
This problem appears (implicitly) in Rothberger's works 
as early as in the 1940's (see, e.g., \cite{ROTH2}).

The terminology in the modern formulation of this problem seems to be 
due to van~Douwen \cite{vD},
where the involved cardinal characteristics of the continuum are named
and studied together with several other cardinal characteristics of the continuum.
All families mentioned below are of infinite sets of natural numbers.
A family $\cF$ is \emph{centered} if the intersection of
each (nonempty) finite subfamily of $\cF$ is infinite.
$\p$ is the minimal size of a centered family which has no \psin{}.
A family $\cF$ is a \emph{tower} if it is linearly quasiordered
by $\as$, and it has no \psin{}.
$\t$ is the minimal size of a tower. 

Clearly, $\p\le\t$.
Thus the Minimal Tower Problem \ref{MTP} is whether $\p=\t$.
$\p$ is starring, implicitly and explicitly, in many constructions carried in the
classical and modern eras of SPM. 
For example, Bell proved that the hypothesis $\p=2^{\aleph_0}$ is equivalent to 
Martin's Axiom for $\sigma$-centered partially ordered sets,
and Fremlin observed that in general, $\p$ is the minimal cardinality
where Martin's Axiom for $\sigma$-centered posets fails \cite{FREMLIN}.

van~Douwen listed in his survey paper \cite{vD} several open problems with regards to the
relationships among various cardinal characteristics of the continuum. All these problems
where solved in the years which followed, except for two really stubborn ones:
Whether $\p<\t$ is consistent (MTP), and whether $\d<\mathfrak{a}$ is consistent.
These questions (especially, the MTP) survived all attempts by almost all mathematicians working
in infinite combinatorics and forcing. The only obtained results were ones asserting
that (most) existing kinds of forcing notions and techniques cannot settle this problem.
One of the problems was that in most extensions by advanced forcing techniques,
the size of the continuum is $\aleph_2$. But it is known that if $\p=\aleph_1$
then $\t=\aleph_1$ too, so $\p=\t$ in these extensions.
Even Shelah has addressed this problem in the past, unsuccessfully. 
A posteriori we can tell that the mathematics was not mature enough then
in order to solve this deep problem.

The major breakthrough came in 2001, when Shelah developed a new method of 
forcing (iteration over a non well-founded set) to show that 
consistently, $\d<\mathfrak{a}$, and thus solve one of the two
remaining problems from van~Douwen's survey paper.
This solution was presented by J\"org Brendle in a series of talks
in a mathematics symposium in honor of Shelah, held in the Ben-Guryon
University. The MTP remained open.

Roughly at the same period, I became interested in topological notions 
related to the MTP. The roots of this interest appeared in my Master's thesis (1997):
Rec\l{}aw showed that the $\gamma$-sets (sets for which every $\w$-cover contains a $\gamma$-cover)
can be thought of as the topological counterpart of the cardinal $\p$. 
Motivated by this, I defined the notion of $\tau$-sets (every $\tau$-cover is a $\gamma$-cover)
so that it becomes the topological counterpart of $\t$, and asked whether these notions coincide
(``No'' if $\p<\t$). Shelah gave a negative solution in ZFC: The Cantor set is a $\tau$-set
(see \cite{tau}). In fact, his argument showed that the Baire space, and therefore any analytic set,
is a $\tau$-set. The solution to this ``approximation problem'' turned out too easy.

Scheepers (personal communication, 1999) suggested that we consider a tighter approximation problem: We know that 
the critical cardinality of $\sone(\Tau,\Gamma)$ (which is stronger than being a $\tau$-set
but is implied by $\sone(\Omega,\Gamma)$, which is the same as being a $\gamma$-set) is
is also $\t$. Thus it is questionable whether this new property is the same as being a 
$\gamma$-set (``No'' if $\p<\t$). This question was also solved in the negative \cite{tautau},
a solution on which I gave a lecture in Shelah's logic seminar at the Hebrew University.
In that lecture (March 2001), I also described a \emph{combinatorial} ``approximation'' of the MTP
which arose from the topological inquiries: Let $\kappa_{\w\tau}$ denote the critical 
cardinality of the property $\binom{\Omega}{\Tau}$ (every $\w$-cover contains a $\tau$-cover).
This cardinal has a simple combinatorial definition (see \emph{Problem of the month} below),
and $\p=\min\{\kappa_{\w\tau},\t\}$.
The question was whether $\p=\kappa_{\w\tau}$. Shelah addressed this problem after the lecture,
and solved it in about 15 minutes (he was not the first one to see this problem!) -- see \cite{ShTb768}.

Later (2003) I started my post-doctoral studies under Shelah's supervision at the Hebrew University.
One of the things I worked on \cite{tautau} 
was a variant of the notion of $\tau$-covers ($\tau^*$-covers) which is closed under
taking de-refinements, and from this a similar problem arose (see \emph{Problem of the Month} below).
When I showed this problem to Shelah he seemed to have had enough of these ``approximation problems''
and told me something like (rough translation from Hebrew): 
``What really should be done is to solve the original problem''.
This conversation was made roughly at the beginning of the current academic year.
Since then a while has passed, and recently I learned that Shelah dedicated the last several
weeks to this problem -- and succeeded.

Shelah's solution is by introducing yet another method of forcing, 
which is a close relative of his solution to the $\d<\mathfrak{a}$ problem
(non well-founded iterations).
Andrzej Roslanowski, who saw the details of the proof,
is about to announce these news in a coming conference.

To our opinion, the importance of the MTP goes much beyond the technique invented in
order to solve it. This problem inspired the imagination of many mathematicians and
led to mathematical inquiries on the forcing method which did not solve it, but turn out very fruitful
in other respects. Moreover, the role of the MTP in topology and the studies of $\tau$-covers
inspired by it have barely begun.

\by{Boaz Tsaban}{tsaban@math.huji.ac.il}

\section{Research announcements}

\subsection{CON$(\p<\t)$}
(See Section \ref{MTP} above.)
\by{Saharon Shelah}{shelah@math.huji.ac.il}

\section{Problem of the month}
Actually the following problem is in pure infinite combinatorics, but
we chose it because: (1) it was motivated by studies of topological diagonalizations
(related to $\tau$-covers), and (2) it is related to the minimal tower problem in a 
straightforward manner.

Let $\inf$ denote the collection of all infinite subsets of $\N$.
Definitions which are missing below can be found in Section \ref{MTP} above.

Let $\kappa_{\w\tau}$ be the minimal cardinality of a centered family $\cF\subseteq\inf$
such that for each $A\in\inf$, the restriction
$\{B\cap A : B\in \cF\}$ is \emph{not} a linearly $\as$-quasiordered subset of $\inf$.
This cardinal was introduced in \cite{tautau}, where it was observed that,
easily, $\p=\min\{\kappa_{\w\tau},\t\}$.
In \cite{ShTb768} it is proved that in fact, $\p=\kappa_{\w\tau}$.
Further study of the notion which yielded that observation
led to the following problem.

\begin{defn}
A family $\cF\subseteq\inf$ is \emph{linearly refinable}
if for each $A\in\cF$ there exists an infinite subset
$\hat A\subseteq A$ such that the family $\hat\cF = \{\hat A : A\in\cF\}$ is
linearly $\as$-quasiordered.
$\p^*$ is the minimal size of a centered family
in $\inf$ which is not linearly refinable.
\end{defn}

Here too, it is easy to see that $\p = \min\{\p^*,\t\}$,
and in \cite{tautau} it is shown that $\p\le\d$.
By Shelah's theorem, it is consistent that $\p^*<\t$.
We therefore have the following problem.
\begin{prob}[\cite{tautau, ShTb768}]
Does $\p=\p^*$?
\end{prob}

Observe that Shelah's solution of the Minimal Tower problem cannot settle this
problem as is, because $\p<\t$ implies $\p=\p^*$.
However, it may well be that this problem can be solved in the positive
by elementary means, or in the negative by standard forcing methods.

\by{Boaz Tsaban}{tsaban@math.huji.ac.il}

\end{document}